\documentclass[12pt]{amsart}

\usepackage{times}
\usepackage[latin1]{inputenc}
\usepackage[T1]{fontenc}
\usepackage{amsmath,amssymb,amsthm}
\usepackage{graphicx}
\usepackage{epsfig}
\usepackage{hyperref}

\usepackage[square,authoryear,sort&compress]{natbib}

\begin{document}

\newtheorem{theorem}{Theorem}[section]
\newtheorem{lemma}[theorem]{Lemma}
\newtheorem{corollary}[theorem]{Corollary}
\newtheorem{proposition}[theorem]{Proposition}
\theoremstyle{definition}
\newtheorem{definition}[theorem]{Definition}
\newtheorem{example}[theorem]{Example}
\newtheorem{claim}[theorem]{Claim}
\newtheorem{xca}[theorem]{Exercise}

\newcommand{\ind}{\mbox{ind}}
\newcommand{\conn}{\mbox{connectivity}}
\newcommand{\Hom}{\mbox{Hom}}

\theoremstyle{remark}
\newtheorem{remark}[theorem]{Remark}
\newcommand{\be}{\begin{equation}}
\newcommand{\ee}{\end{equation}}
\newcommand{\lexmin}{\mbox{lexmin}}
\newcommand{\prob}{\mbox{\bf P}}
\newcommand{\nerve}{\mathcal{N}}
\newcommand{\Star}{\mbox{st}}
\newcommand{\supp}{\mbox{supp}}
\newcommand{\vsupp}{\mbox{vsupp}}
\newcommand{\link}{\mbox{lk}}
\newcommand{\Z}{\mathbb{Z}}
\newcommand{\R}{\mathbb{R}}
\newcommand{\Q}{\mathbb{Q}}
\newcommand{\Homology}{\widetilde{H}}
\newcommand{\bound}{\partial}
\newcommand{\XG}{X(n,p)}
%\theoremstyle{remark}
%\newtheorem{remark}[theorem]{Remark}

%    Blank box placeholder for figures (to avoid requiring any
%    particular graphics capabilities for printing this document).
%\newcommand{\blankbox}[2]{%
%  \parbox{\columnwidth}{\centering
%%    Set fboxsep to 0 so that the actual size of the box will match the
%%    given measurements more closely.
%    \setlength{\fboxsep}{0pt}%
%    \fbox{\raisebox{0pt}[#2]{\hspace{#1}}}%
%  }%
%}

\bibliographystyle{plain}

\title{Topology of random clique complexes}

\author{Matthew Kahle}
\address{Department of Mathematics, University of Washington, Seattle}
\email{mkahle@gmail.com}
\thanks{Supported in part by NSA  \#H98230-05-1-0053 and NSF-VIGRE}

\date{June 11, 2007}

\maketitle

\vspace*{-0.5in}

\begin{abstract}
\noindent In a seminal paper, Erd\H{o}s and R\'enyi identified a
sharp threshold for connectivity of the random graph
$G(n,p)$. In particular, they showed that if $p \gg \log{n}/n$ then $G(n,p)$ is
almost always connected, and if $p \ll \log{n}/n$ then $G(n,p)$ is
almost always disconnected, as $n \rightarrow \infty$.

The {\it clique complex} $X(H)$ of a graph $H$ is the simplicial
complex with all complete subgraphs of $H$ as its faces. In
contrast to the zeroth homology group of $X(H)$, which
measures the number of connected components of $H$, the higher
dimensional homology groups of $X(H)$ do not correspond to monotone graph properties. There are nevertheless higher dimensional
analogues of the Erd\H{o}s-R\'enyi Theorem.

We study here the higher homology groups of $X(G(n,p))$. For $k > 0$ we show the following. If $p=n^\alpha$, with $\alpha < -1/k$ or $\alpha > - 1/(2k+1) $, then the $k$th homology group of $X(G(n,p))$ is almost always vanishing, and if $-1/k < \alpha < -1/(k+1)$, then it is almost always nonvanishing.

We also give estimates for the expected rank of homology, and exhibit explicit nontrivial classes in the nonvanishing regime. These estimates suggest that almost all $d$-dimensional clique complexes have only one nonvanishing dimension of homology, and we cannot rule out the possibility that they are homotopy equivalent to wedges of spheres.
\end{abstract}

\section{Introduction}

A pioneering result in the theory of random graphs is the Erd\H{o}s-R\'{e}nyi theorem on the threshold for connectivity. \cite{Erd1}. This is a primary inspiration for the line of research pursued here, and some of our results may be viewed as generalizations of the Erd\H{o}s-R\'{e}nyi theorem to higher dimensions, so we begin by defining random graphs and stating their result. 

The random graph $G(n,p)$ is defined to be the probability space of all graphs on vertex set $[n]=\{ 1,2, \ldots, n \}$ with each edge inserted independently with probability $p$. Frequently, $p$ is a function of $n$, and one asks whether a typical graph in $G(n,p)$ is likely to have a particular property as $n \rightarrow \infty$. We say that $G(n,p)$ {\it almost always (a.a.)} has property $\mathcal{P}$ if $\mbox{Pr}[G(n,p) \in \mathcal{P}] \rightarrow 1$ as $n \rightarrow \infty$.

\begin{theorem}[{Erd\H{o}s and R\'{e}nyi}] \label{clique_Theorem ER} If
$p=(\log{n}+\omega(n))/n$ and $\omega(n) \rightarrow \infty$ as $n
\rightarrow \infty$ then $G(n,p)$ is almost always connected. If
$\omega(n) \rightarrow - \infty$ then $G(n,p)$ is almost always
disconnected.
\end{theorem}

The number of connected components in a graph is a {\it monotone graph property}. In other words, adding edges to a graph can only decrease the number of components. (As with functions $f:\R \to \R$, we could talk about graph properties either being monotone increasing or decreasing.) Much random graph theory is concerned with monotone graph properties: chromatic number, clique number, subgraph containment, diameter, and so on \cite{Bollo}. 

The {\it clique complex $X(G)$} of a graph $G$ is the simplicial
complex with all complete subgraphs of $G$ as its faces. The $1$-skeleton of $X(G)$ is $G$ itself, so Erd\H{o}s and R\'enyi's result may be interpreted as a statement about homology  $\Homology_0(X(G(n,p)))$ or homotopy  $\pi_0(X(G(n,p))$. (For a brief introduction to the topological notions discussed in this article, please see Section \ref{topology}.)\\

{\bf  Note: To streamline notation, we will abbreviate $X(G(n,p))$ by $\XG$ for the rest of the article.}\\

Our main objects of study are $\Homology_k(\XG)$ and $\pi_k(\XG)$ for each fixed $k>0$.  We find that vanishing of higher homology is not monotone, as homology vanishes for large and small functions $p$, but is nonvanishing for some regime in between. Still, it is possible to make statements which generalize Theorem \ref{clique_Theorem ER}.

%We emphasize that in contrast to the $k=0$ case, the vanishing of higher homology groups in the clique complex is {\it not} a monotone property of the underlying graph. (This is an immeditate corollary of our results, but it is not hard to convince oneself that monotone fails when $k=1$, for example, with small examples.) There are nevertheless analogues of Theorem \ref{clique_Theorem ER}. In particular,

Another way to state our results is to fix the dimension $d$ of the clique complex (by appropriately choosing $p$), rather than looking at a fixed homology group. In this case, we find that the homology of $\XG$ is highly concentrated in its middle dimensions. Asymptotically, a $d$-dimensional random clique complex a.a. has trivial homology above dimension $\lfloor d/2 \rfloor$ and below dimension $\lfloor d/4 \rfloor$. On the other hand, homology is almost always nontrivial in dimension $\lfloor d/2 \rfloor$.

In fact we cannot rule out the possibility that for random clique complexes of bounded dimension, the {\it only} nontrivial homology is in dimension $\lfloor d/2 \rfloor$. We give evidence for this conjecture by estimating the expectation of the rank of homology. It is certainly true that the nontrivial homology in dimension $\lfloor d/2 \rfloor$ accounts for ``almost all'' of the Euler characteristic.

In Section \ref{rsc} we briefly survey other papers concerning topology of random simplicial complexes.  

\section{Topological notions} \label{topology}

The reader who is familiar with reduced homology and homotopy groups of topological spaces may feel free to skip this section. For anyone not so familiar, this will only serve as the very briefest of introductions, and will probably not be sufficient to understand the more technical parts of the proofs, but at the suggestion of an anonymous referee, we are including this section in order to make the article accessible to a wider audience. For anyone who wants to know more, a very nice introduction to algebraic topology is Allan Hatcher's book \cite{Hatcher}.  

The {\it reduced homology groups} $\Homology_i(X,k)$, $i=0,1,2, \ldots$, where $k=\Z$ or some field, are topological invariants associated with a topological space $X$. Very roughly, $\Homology_i(X,k)$ measures the number of $i$-dimensional holes in $X$. Suppose $X$ is a finite simplicial complex of dimension $d$. The most important topological facts for the purposes of this article are the following.

\begin{itemize}

\item $\Homology_i(X,k)$ is a finitely generated abelian group. In the case the $k$ is a field, it is a vector space over $k$.\\

\item The {\it $i$th Betti number} is $\beta_i=\dim{\Homology_i(X,\Q)}$. A classical fact is that if $f_i$ is the number of $i$-dimensional faces of $X$, then the following Euler formula holds. $$f_0 - f_1 +   \cdots + (-1)^d f_d = \beta_0 - \beta_1 + \cdots + (-1)^d \beta_d.$$ Also, it follows directly from the definition of simplicial homology and dimensional considerations that for every $i$, $$-f_{i-1} + f_{i} -f_{i+1} \le \beta_i \le f_{i}.$$

\item $\Homology_0(X,k)=0$ if and only if $X$ is connected.\\

\item $\Homology_i(X,k)=0$ for $i>d$.\\

\item (Universal coefficients for homology) If $\Homology_i(X,\Z)=0$ for $0 \le i \le m$ then $\Homology_i(X,k)=0$ for any coefficients $k$.\\

\end{itemize}

We also briefly discuss the {\it homotopy groups} $\pi_i(X)$. Again, see \cite{Hatcher} for a nice introduction, but the following facts will be more than sufficient to read this article.\\

\begin{itemize}

\item  $\pi_i(X)$ is the set of homotopy classes of maps from the sphere $S^i \mapsto X$. In particular, we say $\pi_0(X)=\{0\}$ if and only if $X$ is path connected.\\

\item (Hurewicz Theorem) If $\pi_i(X)=\{0\}$ for $i \le n$,  (in which case we say that $X$ is {\it $n$-connected}) then $\Homology_i(X,\Z)=0$, $i=0, 1, \ldots, n$.\\

\item $\Homology_1(X,\Z)$ is the abelianization of the fundamental group $\pi_1(X)$.\\

\end{itemize}

Reduced homology groups and homotopy groups are topological invariants, meaning that if two spaces are homeomorphic then their associated homology and homotopy groups are isomorphic. A stronger statement, also true, is that they are homotopy invariants, meaning that the same holds even if the spaces are only homotopy equivalent.

\section{Statement of results}

We discuss which groups $\widetilde{H}_i(\XG,\Z)$ are nontrivial, then estimate Betti numbers. For comparison with the results, note that $$\dim(\XG) \approx -2\log{n}/\log{p}.$$ For example, since $\dim(X(H)) \ge k$ if and only if $H$ contains $(k+1)$-cliques,  standard random graph techniques for subgraph containment \cite{Bollo} give that if $p=n^{\alpha}$ with $\alpha < -2/k$ then a.a. $\dim(\XG)< k,$ and if $\alpha > -2/k$ then a.a. $\dim(\XG) \ge k$.

We first show that if $p$ is large enough then homology vanishes.
A topological space $X$ is said to be {\it $k$-connected} if every map from a sphere $\mathbb{S}^{i}
\rightarrow X$ extends to a map from the ball $\mathbb{B}^{i+1}
\rightarrow X$ for $i =0, 1, \ldots, k$. Equivalently, $X$ is $k$-connected if $\pi_i(X)=0$ for $i \le k$, and in particular $0$-connected is equivalent to path connected. This implies, by the Hurewicz Theorem \cite{Hatcher}, that $\widetilde{H}_i(X,\Z)=0$ for $i \le k$.

The following is implicit in \cite{Mesh}, although Meshulam's result was more general and stated for homology instead of homotopy groups. We prove the homotopy statement here for the sake of completeness, although the argument is similar to Meshulam's.

\begin{theorem}[Meshulam] \label{clique_Theorem M} If every $2k+2$ vertices of a graph $H$ have a common neighbor then $X(H)$ is $k$-connected.
\end{theorem}

In the case that $H$ is a random graph, this can be improved. For example, in the case $k=0$, Erd\H{o}s and R\'enyi's theorem gives that the
threshold for connectivity is the same as the threshold for every vertex having at least one neighbor. The threshold for every set of $l$ vertices having a neighbor is given by the following.

\begin{theorem}\label{clique_Theorem 1} If $p=\left( \frac{l \log{n} + \omega(n)}{n}
\right)^{1/l}$ and $\omega(n) \rightarrow \infty$ then a.a., every
$l$ vertices of $G(n,p)$ have a common neighbor.
\end{theorem}

Together with Meshulam's result, we immediately have the following.

\begin{corollary}\label{Corollary 1} If $p=\left( \frac{(2k+2) \log{n} + \omega(n)}{n}
\right)^{1/(2k+2)}$ and $\omega(n) \rightarrow \infty$ then a.a.
$\XG$ is $k$-connected.
\end{corollary}

Corollary \ref{Corollary 1} can be improved and we do so with Theorem \ref{clique_Theorem 2}. Note that when $k=0$, this specializes to one direction
of the Erd\H{o}s-Renyi theorem.

\begin{theorem}\label{clique_Theorem 2} If $p=\left( \frac{(2k+1) \log{n} + \omega(n)}{n}
\right)^{1/(2k+1)}$ and $\omega(n) \rightarrow \infty$ then a.a.
$\XG$ is $k$-connected.
\end{theorem}

As a consequence, we have a statement about vanishing of homology. In a different regime, we can make statements about nonvanishing homology by exhibiting nontrivial classes explicitly.

\begin{theorem} \label{clique_Theorem 4} If $p^{k+1}n \rightarrow 0$
and $p^{k}n \rightarrow \infty$ as $n\rightarrow \infty$ then
$\XG$ a.a. retracts onto a sphere $\mathbb{S}^{k}$. Hence
$\widetilde{H}_k(\XG,\Z)$ a.a. has a $\Z$ summand.
\end{theorem}

Theorem \ref{clique_Theorem 2} gives a statement that if $p$ is large enough then homology vanishes. The same must be true when $p$ is small enough, simply by dimensional considerations. But this kind of coarse argument will only give that $\alpha < -2/k$ then $\Homology_k(\XG,\Z)=0$. By Theorem \ref{clique_Theorem 4}, the following is best possible.

\begin{theorem} \label{vanish} If $p=n^{\alpha}$
with $\alpha < -1/k$ then $\Homology_k(\XG,\Z)=0$
almost always.
\end{theorem}

By Theorems \ref{clique_Theorem 2}, \ref{clique_Theorem 4}, and \ref{vanish}, we have the following.

\begin{corollary}[Vanishing and nonvanishing of homology] \label{homo} If $p=n^{\alpha}$ then
\begin{enumerate}
\item if $\alpha < -1/k$ or $\alpha > -1/(2k+1)$ then a.a. $\Homology_k(X(G(n,p),\Z)=0$,
\item and if $-1/k < \alpha < -1/(k+1)$ then a.a. $\Homology_k(\XG,\Z) \neq 0$.
\end{enumerate}
\end{corollary}

So rather than monotonicity, we have a kind of unimodality (in terms of $p$) for each fixed homology group as $n \to \infty$.

Corollary \ref{homo} does not address the case when $- 1/(k+1) < \alpha < -1/(2k+1)$. We believe that Theorem \ref{clique_Theorem 2} can probably be improved to say that if $p=n^{\alpha}$ with $\alpha > -1/(k+1)$ then a.a. $\Homology_k(\XG,\Z)=0$.

To give evidence for this conjecture, we estimate the expected rank of homology, and show that it passes through phase transitions near $\alpha=-1/k$ and $-1/(k+1)$. Let $f_k$ denote
the number of $k$-dimensional faces of $\XG$ and
$\beta_k$ its $k$th Betti number. That is, let $$\beta_k=\dim{\Homology_k(\XG,\Q)},$$ although our result holds for coefficients in any field. By the definition of simplicial homology and dimensional considerations, $\beta_k \le f_k$.

We show that given the hypothesis of Theorem \ref{clique_Theorem 4}, $f_k$ is actually a good approximation for $\beta_k$,
but for $p$ outside of this range, $\beta_k$ is much smaller. 
We write $X \sim Y$ almost always if for every $\epsilon > 0$, as $n \to \infty$,
$$\prob((1-\epsilon) \le Y/X \le (1+\epsilon)) \to 1.$$

\begin{theorem} \label{clique_Theorem 5} If $p^{k+1}n \rightarrow 0$
and $p^{k}n \rightarrow \infty$ then
$E(\beta_k)  / E(f_k)  \to 1$. Moreover $\beta_k \sim E[\beta_k]$ and $f_k \sim E[f_k]$ a.a., so  $\beta_k \sim f_k$ a.a.
\end{theorem}

Finally, we apply discrete Morse theory to show that $E[\beta_k]/E[f_k]$ passes through phase transitions at $p=n^{-1/(k+1)}$ and $p=n^{-1/k}$.

\begin{theorem} \label{clique_Theorem 6} If $p^{k+1}n \rightarrow \infty$
 or $p^{k}n \rightarrow 0$ then $E(\beta_k)/E(f_k) \to 0$.
\end{theorem}

(Note that even the second case of Theorem \ref{clique_Theorem 6} is not necessarily implied by Theorem \ref{vanish}, since the statement that a random variable is a.a. zero implies nothing about its expectation. Also, $p^{k}n \to 0$ is a slightly weaker hypothesis than $p=n^{\alpha}$ with $\alpha < -1/k$.)

As a corollary to Theorems \ref{clique_Theorem 5} and \ref{clique_Theorem 6} we have the following.

\begin{corollary}[Betti numbers] If $p=n^{\alpha}$ then for any $\epsilon > 0$,
\begin{enumerate}
\item if $\alpha < -1/k$ or $\alpha > -1/(k+1)$ then a.a. $0 \le \beta_k / f_k < \epsilon$,
\item if $-1/k < \alpha < -1/(k+1)$ then a.a. $ 1 - \epsilon < \beta_k / f_k \le 1 $.
\end{enumerate}
\end{corollary}

If Theorem \ref{clique_Theorem 2} can be improved to say that if $\alpha > -1/(k+1)$ then $\XG$ is a.a. $k$-connected, then the upshot is that a.a. $d$-dimensional clique complexes have only one nonvanishing dimension of homology. This might be a bit surprising, since it does not depend on $d$. In a sense this almost determines the homotopy type. We discuss this more in Section \ref{open}.

In the next several sections we prove the results. Theorems \ref{clique_Theorem M}, \ref{clique_Theorem 1} and \ref{clique_Theorem 2} are proved in Section \ref{connectivity}, Theorem \ref{vanish} in Section \ref{vanishing}, Theorem \ref{clique_Theorem 4} in Section \ref{spherical}, Theorem \ref{clique_Theorem 5} and Theorem \ref{clique_Theorem 6} in Section \ref{Betti 1}.

\section{Connectivity} \label{connectivity}

We use the following Nerve Theorem of Bj\"{o}rner \cite{Bjorner} throughout this section. The {\it nerve} of a family of nonempty sets $(\Delta_i)_{i \in I}$ is the simplicial complex $\nerve((\Delta_i)_{i \in I})$, defined on the vertex set $I$ by the rule that $\sigma \in \nerve(\Delta_i) $ if and only if $\cap_{i \in \sigma} \Delta_i \neq \emptyset$. Note that the nerve depends on the whole family, but for brevity's sake we denote it by $\nerve(\Delta_i)$ rather than $\nerve((\Delta_i)_{i \in I})$.

\begin{theorem}[Bj\"{o}rner]\label{nerve}
Let $\Delta$ be a simplicial complex, and $(\Delta_i)_{i \in I}$ a family of subcomplexes such that $\Delta = \cup_{i \in I} \Delta_i$. Suppose that every nonempty finite intersection $\Delta_{i_1} \cap  \Delta_{i_2} \cap \ldots \cap \Delta{i_t}$ is $(k-t+1)$-connected, $t \ge 1$. Then $\Delta$ is $k$-connected if and only if $\nerve(\Delta_i)$ is $k$-connected.
\end{theorem}

\begin{proof}[Proof of Theorem \ref{clique_Theorem M}]

We show that if every $2k+2$ vertices of a graph $H$ have a neighbor then $X(H)$ is a.a. $k$-connected.
Proceed by induction on $k$. The claim holds when $k=0$, i.e. if every pair of vertices of a graph have a common neighbor, then the graph is certainly connected. So suppose the claim holds for $k=0, \ldots, i-1$ where $i \ge 1$. Further suppose that $H$ is a graph such that every set of $2i+2$ vertices has some neighbor. We wish to show that $X(H)$ is $i$-connected.

Define the {\it star} of a vertex $v$ in a simplicial complex $\Delta$ to be the subcomplex $\Star_{\Delta}(v)$ of all faces in $\Delta$ containing $v$. Clearly we have $\Delta=\cup_{v \in \Delta} \Star_{\Delta}(v)$. So to apply Theorem \ref{nerve} we must check that each vertex star is itself $i$-connected, and that every $t$-wise intersection is $(i-t+1)$-connected for $t=2,\ldots,i+1$.

Each star is a cone, hence contractible and in particular $i$-connected. Each $t$-wise intersection of vertex stars is a clique complex in which every $2i+2-t$ vertices share a neighbor, hence by induction is $i-\lfloor t/2 \rfloor$-connected. Since $t \ge 2$, $i- \lfloor t/2 \rfloor \ge i-t+1$, so the claim follows provided that $\nerve(\Star_{X(H)}(v))$ is also $i$-connected. This is clear though; since every $2i+2$ neighbors have a neighbor, the intersection of every $2i+2$ vertex stars is nonempty. So the $(2i+1)$-dimensional skeleton of $\nerve(\Star_{\Delta}(v))$ is complete, and then $\nerve(\Star_{\Delta}(v))$  is $2i$-connected.

\end{proof}

\begin{proof}[Proof of Theorem \ref{clique_Theorem 1}]
We claim that if
$p=\left( \frac{l \log{n} + \omega(n)}{n} \right)^{1/l}$ and
$\omega(n) \rightarrow \infty$ then a.a. every $l$ vertices of
$G(n,p)$ have a neighbor. The expected number of $l$-tuples of vertices
in $G(n,p)$ with no neighbor is

\begin{eqnarray*}
&&{n \choose l}(1-p^l)^{n-l}\\
&\le& {n \choose l} e^{-p^l(n-l)}\\
&=&{n \choose l} e^{-\frac{l \log{n} + \omega(n)}{n} (n-l)}\\
&=&{n \choose l} n^{-l} e^{-\omega(n)(n-l)/n}\\
&\le &  e^{-\omega(n)(1-l/n)}\\
&=&o(1),\\
\end{eqnarray*}

\noindent since $\omega(n) \rightarrow \infty$. This proves Theorem \ref{clique_Theorem 1}.

\end{proof}

For a graph $H$ and any subset of vertices $U \subseteq V(H)$, define $$S(U):=\bigcap_{v \in U} \Star_{X(H)}(v).$$

\begin{lemma}\label{conn lemma} Let $k \ge 1$ and suppose $H$ be any graph such that every $2k+1$ vertices share a neighbor, and for every set of $2k$ vertices $U \subseteq H$, $S(U)$ is connected. Then $X(H)$ is $k$-connected.
\end{lemma}

\begin{proof}[Proof of Lemma \ref{conn lemma}]

As in the proof of Theorem \ref{clique_Theorem 1}, cover $X(H)$ by its vertex stars $\Star(v)$ and apply Theorem \ref{nerve}.  The nerve $\nerve(\Star(v))$ is $k$-connected since every $2k+1$ vertices sharing a neighbor implies that its $2k$-skeleton is complete, so it is in fact $(2k-1)$-connected. Then to check that $X(H)$ is $k$-connected, it suffices to check that every $t$-wise intersection of vertex stars is $(k-t+1)$-connected, $2 \le t \le k+1$.  We show something slightly stronger,  that if $0 \le j < k $ and $i \le 2k-2j$, then every $i$-wise intersection of vertex stars is $j$-connected.

The case $j=0$ is clear: if $|U|=2k$ then $S(U)$ is connected by assumption, and if  $|U| < 2k$ then $S(U)$ is still connected, since every pair of vertices in $S(U)$ shares a neighbor. Let $j=1$. The claim is that if $i \le 2k-2$ and $|U|=i$ then $S(U)$ is $1$-connected. Cover $S(U)$ by vertex stars $\Star_{S(U)}(v)$, $v \in S(U)$ and again apply Theorem \ref{nerve}.  We only need to check that every intersection $\Star_{S(U)}(v) \cap \Star_{S(U)}(v)$ is connected, but this is clear since $$\Star_{S(U)}(u) \cap \Star_{S(U)}(v)=S(U\cup\{u,v\})$$ is the intersection of $i+2 \le 2k$ vertex stars, connected by assumption.

Similarly, let $j=2$, $i \le 2k-4$, and $|U|=i$. Then to show that $S(U)$ is $2$-connected, cover by vertex stars $\Star_{S(U)}(v)$. Each $3$-wise intersection of vertex stars $$\Star_{S(U)}(u) \cap \Star_{S(U)}(v) \cap  \Star_{S(U)}(w) = S(U \cup \{u,v,w\})$$ is the intersection of at most $i+3 \le 2k-1$ vertex stars, connected by assumption. Each $2$-wise intersection of vertex stars in $S(U)$ is the intersection of at most $i+2 \le 2k-2$ vertex stars in $X(H)$, $1$-connected by the above. Again applying Theorem \ref{nerve}, we have that $S(U)$ is $2$-connected as desired.

Proceeding in this way, the lemma follows by induction on $j$.

\end{proof}

\begin{proof}[Proof of Theorem \ref{clique_Theorem 2}]

The remainder of this section is a proof that if $$p=\left( \frac{(2k+1)\log{n} + \omega(n)}{n} \right)^{1/(2k+1)}$$ and $\omega(n)
\rightarrow \infty$ then a.a. $\XG$ is $k$-connected. Our argument is inspired by a proof of Theorem \ref{clique_Theorem ER} in \cite{Bollo}. Since $k=0$ is Theorem \ref{clique_Theorem ER} we assume that $k \ge 1$,  Observe that for any graph $H$ and vertex subset $U \subseteq V(H)$, $S(U)$ is the clique complex of a subgraph of $H$. Moreover, for any vertex $v \in S(U)$, $\Star_{S(U)}(v)=S(U \cup \{ v \})$. We use these facts repeatedly.

By Lemma \ref{conn lemma} and Theorem \ref{clique_Theorem 1}, we need only check that a.a. the intersection of every $2k$ vertex stars in $\XG$ is connected.  It is convenient to instead check that the intersection of every $2k$ vertex links is connected. For a vertex $v$ in a simplicial complex $\Delta$ define the {\it link} of $v$ in $\Delta$ by $$\link_{\Delta}(v):=\{ \sigma | v \notin \sigma \mbox{ and } \{v\} \cup \sigma \in \Delta \},$$ and for any vertex set $U$ denote $$L(U):=\bigcap_{v \in U}\link(v).$$

Suppose $H$ is as in the hypothesis of Lemma \ref{conn lemma} and $|U|=2k$. If $L(U)$ is connected, then $S(U)$ is connected also, as follows. If $S(U)-L(U)=\emptyset$ we are done, so suppose $x \in S(U)-L(U)$. Clearly $x \in U$. $L(U)$ is connected by assumption, and in particular nonempty, so let $v \in L(U)$. For $u \in U-\{x\}$, $v \sim u$ and $x \sim u$. So $\{ u,v \}$, $\{u,x\}$, and $\{v,x\}$ are all edges in $H$, and $\{u,v,x\}$ is a face in $X(H)$, and $\{ v,x \} \in \Star_{X(H)}(u)$. So $\{ v,x \} \in S(U)$ and $x$ is connected to $L(U)$. This holds for every $x \in S(U)-L(U)$, so $S(U)$ is connected.

Now we check that if $$p=\left( \frac{(2k+1)\log{n} + \omega(n)}{n} \right)^{1/(2k+1)},$$ then a.a., for every subset $U \subseteq [n]$ with $|U|=2k$, $L(U)$ is connected. It suffices to consider the $1$-dimensional skeleton $L(U)^{(1)}$, which is a random graph with independent edges. However the number of vertices in the graph is not constant but a distribution, and there are ${n \choose 2k}$ such graphs, where edges in one are not necessarily independent of edges in another. However, the edges within each graph are still independent, and we may still apply linearity of expectation to show that the probability that at least one of these graphs is not connected goes to $0$. 

Let $U\subseteq [n]$ be any vertex set of cardinality $2k$. The number of vertices $X$ in $L(U)$ is not constant, but it is tightly concentrated. $X$ is the sum of $n-2k$ independent indicator random variables, each with probability $p^{2k}$. So we have an the following estimate for the mean of $X$. $$\mu = E[X] \sim p^{2k}n$$ since $k$ is constant. It is convenient to assume that $p = o(1)$. A similar argument works for dense random graphs.

Standard large deviation bounds \cite{Alon} give that $$\prob (|X-\mu| > \epsilon \mu ) < e^{-c_{\epsilon}\mu}$$ for some constant $c_{\epsilon}>0$ depending only on $\epsilon$. We set $\epsilon=1/100$ and write $c=c_\epsilon$. Then

\begin{eqnarray*}
e^{-c \mu} & \le & e^{-cp^{2k}n} \\
& = & e^{-cp^{-1}p^{2k+1}n} \\
& \le & e^{-cp^{-1}(2k+1)\log{n}} \\
& \le & n^{-c(2k+1)p^{-1}}\\
& \le & n^{-c(2k+1)\omega(n)}\\
\end{eqnarray*} where $\omega(n) \to \infty$. So, applying a union bound, the total probability that for any set $U$, $|X-p^{2k}n|>(1/100)p^{2k}n$ is no more than $${n \choose 2k} n^{-c(2k+1)\omega(n)} = o(1).$$ We have shown that a.a., $0.99p^{2k}n < X < 1.01p^{2k}n$ holds for every $U$, so we assume this for the remainder of the proof. Note that $X \to \infty$ by our assumption on $p$.

Let $\prob_i$ denote the probability that there are components of order $i$ in $L(U)$ for at least one $2k$-subset $U$. $\prob_1 =o(1)$ by Theorem \ref{clique_Theorem 1}. Next we bound $\prob_2$. There are ${n \choose 2k}$ choices for $U$, then conditioned on that choice of $U$, let $X$ denote the number of vertices in $L(U)$, as above. Given that $u,v \in L(U)$, the probability that $\{u,v\}$ spans a component of order $2$ in $L(U)$ is $p(1-p)^{2(X-2)}$. There are ${X \choose 2}$ choices for $\{u,v\}$ so by our assumptions on $X$,

\begin{eqnarray*}
\prob_2 & \le & {n \choose 2k}{X \choose 2} p(1-p)^{2(X-2)}\\
& \le & n^{2k} {\lceil 1.01p^{2k}n \rceil \choose 2} p e^{-2p(X-2)}\\
& \le & n^{2k} p^{4k}n^2 p e^{-2p(X-2)}\\
& \le & n^{2k+2} p^{4k+1}e^{-2pX(1-o(1))}\\
& \le & n^{2k+2} p^{4k+1}e^{-1.98p^{2k+1}n(1-o(1))}\\
& \le & n^{2k+2} p^{4k+1}e^{-1.98(2k+1)\log{n}(1-o(1))}\\
& \le & n^{2k+2} p^{4k+1}n^{-1.98(2k+1)(1-o(1))}\\
& = & o(n^{-1}),\\
\end{eqnarray*} since $k \ge 1$.

For any $U$ and subset of $i$ vertices $S \subseteq L(U)$, for $S$ to span a connected component of order $i$, it must at least contain a spanning tree. It is well known that the number of spanning trees on $i$ vertices is $i^{i-2}$ \cite{StanleyEC2}. The probability that all $i-1$ edges in any particular tree appear is $p^{i-1}$, by independence. We first bound $P_i$ from above, assuming $3 \le i \le 100$. Since $X \to \infty$,

\begin{eqnarray*}
\prob_i & \le & {n \choose 2k}{X \choose i}i^{i-2} p^{i-1}(1-p)^{i(X-i)}\\
&\le & n^{2k} \frac{X^i}{i!} i^{i-2} p^{i-1} e^{-ipX(1-o(1))}\\
&\le & c_i n^{2k} X^i p^{i-1} e^{-i(0.99p^{2k+1}n(1-o(1))}\\
&\le & c_i n^{2k} (1.01p^{2k}n)^i p^{i-1} e^{-0.99i(2k+1)\log{n}(1-o(1))}\\
%&= & n^{2k+i} 1.01^i p^{2ki+i-1} e^{-piX(1-o(1))}\\
%&\le & n^{2k+i} 1.01^i p^{2ki+i-1} e^{-piX(1-o(1))}\\
&= & c_i  \exp [(2k+i)\log{n} + i\log{1.01} + (2ki+i-1)\log{p}\\
& & -0.99(1-o(1))i(2k+1)\log{n}]\\
&\le & c_i  \exp[(2k+i -0.99i(2k+1)+o(1))\log{n}]\\
&\le & c_i \exp[(2k+0.01i -1.98ik+o(1))\log{n}]\\
&\le & c_i \exp[(2k+0.01i -ik-0.98ik+o(1))\log{n}]\\
&\le & c_i \exp[(-k-0.97i + o(1))\log{n}],\\
\end{eqnarray*} where $c_i=i^{i-2}/i!$ is a constant that only depends on $i$. (The last line holds because $i \ge 3$ and $k \ge 1$.) So for large enough $n$,

\begin{eqnarray*}
\prob_i & \le & c_i \exp[(-k-0.97i + o(1))\log{n}].\\
&\le & c_i n^{-k/2-.97i}\\
\end{eqnarray*}
and $$\sum_{i=3}^{100}\prob_i \le \sum_{i=3}^{100} c_i n^{-k/2-.97i} = o(n^{-3}).$$

Now suppose $100 < i \le \lfloor 0.6p^{2k}n \rfloor$. Here we need to be a bit more careful in our treatment of the $i^{i-2}/i!$ factor. Stirling's formula gives that $i^{i-2}/i! \le e^i$ though, and this will be good enough. We have

\begin{eqnarray*}
\prob_i & \le & {n \choose 2k}{X \choose i}i^{i-2} p^{i-1}(1-p)^{i(X-i)}\\
&\le & n^{2k} \frac{X^i}{i!} i^{i-2} p^{i-1} e^{-p(0.4iX)}\\
&\le & n^{2k} X^i e^i p^{i-1} e^{-0.4ipX}\\
&\le & n^{2k} (1.01p^{2k}n)^i e^i p^{i-1} e^{-0.4i(0.99p^{2k+1}n)}\\
%&= & n^{2k+i} 1.01^i p^{2ki+i-1} e^{-(0.396ip^{2k+1}n)}\\
%&\le & n^{2k+i} 1.01^i p^{2ki+i-1} e^{-0.396i(2k+1)\log{n}}\\
&= & \exp [(2k+i)\log{n} + i(1+\log{1.01}) + (2ki+i-1)\log{p}\\
& & -0.396i(2k+1)\log{n}] \\
&\le & \exp[(2k+(0.604+o(1))i-0.792ik)\log{n}]. \\
\end{eqnarray*}
%since $\log{p} \le 0$ and $\log{1.01}/\log{n}=o(1)$.

Then by assumption that $k \ge 1$, $i > 100$, and for large enough $n$,

\begin{eqnarray*}
\prob_i & \le & \exp[(2k+0.605i-0.792ik)\log{n}]\\
& = & \exp[(2k+0.605i-0.092ik-0.7ik)\log{n}]\\
& \le & \exp[(2k+0.605i-9.2k-0.7i)\log{n}]\\
& = & \exp[(-7.2k-.095i)\log{n}]\\
&=&n^{-7.2k-.095i},
\end{eqnarray*}

and $$\sum_{i=101}^{\lfloor 0.6p^{2k}n \rfloor} \prob_i \le \sum_{i=101}^{\infty} n^{-7.2k-.095i} = o(n^{-15}).$$

Putting it all together, a.a. each $L(U)$ is of order $X$, $0.99p^{2k}n<X<1.01p^{2k}n$, and there are no components of order $i$, $1 \le i \le \lfloor 0.6 p^{2k}n \rfloor $ in any of the $L(U)$. We conclude that each $L(U)$ is connected, as desired.

\end{proof}

\section{Vanishing homology} \label{vanishing}

We show if $p=n^{\alpha}$ with $\alpha < -1/k$ then $\Homology_k(X(G(n,p),\Z)=0$ almost always. 
In this section, we assume the reader is familiar with simplicial homology \cite{Hatcher}.
For a $k$-chain $C$, the
{\it support}, $\supp(C)$, is the union of $k$-faces 
in $C$ with nonzero coefficients. Similarly,
the {\it vertex support}, $\vsupp(C)$, is the underlying vertex set of the support. 

A pure $k$-dimensional subcomplex $\Delta$ is said to be {\it strongly connected} 
if every pair
of $k$-faces $\sigma,\tau \in \Delta^d$ can be connected by a sequence of 
facets $\sigma = \sigma_o, \sigma_1, \sigma_2, \ldots \sigma_j = \tau$
such that $\dim(\sigma_i \cap \sigma_{i+1}))=d-1$ for $0 \le i \le n-1$. 
Every $k$-cycle is a $\Z$-linear combination of $k$-cycles
with strongly connected support. We show first that all strongly connected 
subcomplexes are supported
on a bounded number of vertices, and then that all small cycles are boundaries. 

\begin{lemma} \label{strong}
Let $\alpha < -1/k$ and
$0 < 1/N < -1/k -\alpha$. Then there are a.a. no strongly 
connected pure $k$-dimensional
 subcomplexes of $\XG$ with vertex support of more than $N+k+1$ vertices.
\end{lemma}

\begin{proof}[Proof of Lemma \ref{strong}]
The vertices in the support of a strongly connected subcomplex can be ordered 
$v_1, v_2, \ldots v_n$ such that $\{ v_1,\ldots, v_{k+1} \}$ spans a $k$-face 
and $v_i$ is connected to at least $k$ vertices $v_j$ with $j<i$. One way to see 
this is to order the $k$-faces $f_1, f_2, f_3, \ldots$, so that each has $(k-1)$-dimensional 
intersection with the union of the previous faces. That this is possible is guaranteed 
by the assumption of strongly connected. Then let this ordering induce an ordering 
on vertices, since at most one new vertex gets added at a time
in the sequence $f_1, f_1 \cup f_2, f_1 \cup f_2 \cup f_3 \ldots$

Suppose $\Delta$ has $N+k+1$ vertices. There are at least 
${k+1 \choose 2} + Nk$ edges in $\Delta$ by the above.
If the underlying graph of $\Delta$ is not a subgraph of $G(n,p)$
then $\Delta$ is not a subcomplex. Choose $\epsilon$ and $N$ such that $1/N < \epsilon < -\alpha -1/k$.
We apply a union bound on the total probability that there are 
any subcomplexes isomorphic to $\Delta$ in $\XG$. We have $p = n^{\alpha} < n^{-(1/k+\epsilon)}$ and $k < \epsilon N k$ by assumption, so

\begin{eqnarray*}
\prob(\exists \mbox{ subcomplex}) & \le &(N+k+1)! {n \choose N+k+1} p^{ {k+1 \choose 2} + Nk} \\
& \le & (N+k+1)! {n \choose N+k+1}  n^{-(1/k+\epsilon) ({k+1 \choose 2} + Nk)} \\
& \le & n^{N+k+1}  n^{-(1/k)( {k+1 \choose 2} + Nk)} n^{-\epsilon({k+1 \choose 2} + Nk)}  \\
& \le & n^{1-(k+1)/2 - \epsilon {k+1 \choose 2}} \\
& = & O(n^{-\epsilon}). \\
\end{eqnarray*} This last line holds since $k \ge 1$. There are only finitely many isomorphism 
types of strongly connected $k$-dimensional complexes $\Delta$ on $N+k+1$ vertices, and a.a.
none of them are subcomplexes of $\XG$ by repeating this argument for each of them. 
There are also no such subcomplexes on more than $N+k+1$ vertices, since each of these 
contains a strongly connected subcomplex on exactly $N+k+1$ vertices (e.g., by the 
ordering described above).

\end{proof}

Then homology is generated by cycles supported on small vertex sets.
Let $\gamma$ be a nontrivial $k$-cycle in a simplicial complex $\Delta$,
with minimal vertex support, and write it is a linear combination of faces
$$\gamma=\sum_{f \in \supp(\gamma)} \lambda_{f} f,$$ with $\lambda_f 
\in \Z$. For the remainder of this section we restrict our attention 
to the full induced subcomplex
of $\Delta$ on $\vsupp(\gamma)$.
Clearly $\gamma$ is still a nontrivial cycle in this subcomplex.
For $v \in \vsupp(\gamma)$ define the $k$-chain $$\gamma \cap \Star(v) := \sum_{f \in \Star(v)} \lambda_{f} f,$$
and the $(k-1)$-chain $$\gamma \cap \link(v) := \sum_{f \in \Star(v)} \lambda_{f} (f-\{v\}).$$

Order the vertices with $v$ last and let this induce an orientation
on every face. We observe that $$\gamma \cap \link(v)=\bound(\gamma \cap \Star(v)),$$ and since $\bound \circ \bound =0$
this gives that $\gamma \cap \link(v)$ is a $(k-1)$-cycle.

\begin{lemma} \label{nontriv} With notation as above, $\gamma{ }\cap{ }\link(v)$ 
is a nontrivial $(k-1)$-cycle in $\link(v)$.
\end{lemma}
\begin{proof}[Proof of Lemma \ref{nontriv}] We need only check that $\gamma{ }\cap{ }\link(v)$ is not a boundary.
Suppose by way of contradiction that $\bound(\beta)=\gamma{ }\cap{ }\link(v)$ for some $k$-chain $\beta$ with
$\supp(\beta) \subseteq \link(v)$. In particular $v \notin \vsupp(\beta)$. Write $$\beta=\sum_{f \in \supp(B)} 
\mu_f f$$ with $\mu_f \in \Z$ and
define the $(k+1)$-chain
$$\beta * \{v \} := \sum_{f \in \supp(\beta)} \mu_f (f \cup \{ v \}).$$ Then 
$$\bound(\beta * \{v \})= \gamma\cap\Star(v) + (-1)^{k+2} \beta .$$ 
So $$\gamma':=\left( \gamma - \gamma\cap\Star(v)\right)+ (-1)^{k+3} \beta $$ is a $k$-cycle
homologous to $\gamma$, but
with $\vsupp(\gamma') \subseteq \vsupp(\gamma)-\{v\}$, contradicting 
that $\gamma$ has minimal vertex support.
\end{proof}

\begin{lemma} \label{octa}  Let $H$ be a graph and $X(H)$ its clique complex.
Suppose $\gamma$ is a nontrivial $k$-cycle in $X(H)$. 
Then $|\vsupp(\gamma)| \ge 2k+2$.
\end{lemma}
\begin{proof}[Proof of Lemma \ref{octa}] Proceed by induction on $k$. The claim is clear when $k=0$.
Suppose then that $|\vsupp(\gamma)| \le 2k+1$, and $v \in \vsupp(\gamma)$.
By Lemma \ref{nontriv}, $\gamma \cap \link(v)$ is a nontrivial cycle. By
the induction hypothesis, $| \vsupp(\gamma\cap\link(v)) | \ge 2k$, so we must have 
equalities $| \vsupp(\gamma) | =2k+1$ and 
$| \vsupp(\gamma\cap\link(v)) | = 2k$. Repeating this argument gives that 
every vertex in $\vsupp(\gamma)$ 
has degree $2k$, so $\vsupp(\gamma)$ spans a clique in $H$. 
But then $\vsupp(\gamma)$ spans a $2k$-dimensional face
in $X(H)$, a contradiction to $\gamma$ nontrivial.
\end{proof}

\begin{proof}[Proof of Theorem \ref{vanish}]
Any nontrivial $k$-cycle with
minimal vertex support must have minimum vertex degree
at least $2k$ in its supporting subgraph, since each vertex 
link is a nontrivial $(k-1)$-cycle by
Lemma \ref{nontriv}, hence contains at least $2(k-1)+2=2k$ vertices by
Lemma \ref{octa}. (We discuss nontrivial $k$-cycles $S^k$ with $| \vsupp(S^k) | =2k+2$
in
Section \ref{sphere}.) 

Let $H$ be any fixed graph with minimal vertex degree $2k$. Let $m=|V(H)|$, and then $|E(H)| \ge m(2k)/2=mk$.
Then if $\alpha < -1/k$ and $p=n^{\alpha}$, $H$ is a.a. not a subgraph of $G(n,p)$. We check this with a union bound. The probability that $H$ is a subgraph is at most

\begin{eqnarray*}
& & m! {n \choose m} p^{mk} \\
& \le & n^m n^{\alpha mk} \\
& = & o(1), \\
\end{eqnarray*} since $\alpha k < -1$. There are only finitely many isomorphism types of graphs
of minimal degree $2k$ on $m=N+k$ vertices. Each has at least $km$ edges.
Applying this argument to each of them, 
we conclude
$\XG$ a.a. has no vertex minimal nontrivial $k$-cycles, 
so a.a. $\Homology_k(X(G(n,p),\Z)=0$.
\end{proof}

\section{Spherical retracts} \label{spherical}
\label{sphere}

We prove Theorem \ref{clique_Theorem 4}, that if
$p^{k+1}n \rightarrow 0$ and $p^{k}n \rightarrow \infty$ as
$n\rightarrow \infty$ then $\XG$ a.a. retracts onto a sphere
$\mathbb{S}^{k}$.

Let $S^d$ denote the $d$-dimensional octahedral sphere (i.e. the $d$-fold repeated join of two isolated points), and $(S^d)^{(1)}$
its $1$-skeleton. An alternate description of $(S^d)^{(1)}$ as a
graph is 
$$V((S^d)^{(1)})=\{u_1,u_2, \ldots, u_{d+1}\} \cup \{v_1,v_2, \ldots, v_{d+1} \}$$ and 
$$E((S^d)^{(1)})=- \{ \{u_i,v_j\} \mid i = j \}$$ where the $-$ denotes complement in the
set of all possible edges. Hence $(S^k)^{(1)}$ has $2(k+1)$
vertices and ${2(k+1) \choose 2}-(k+1)$ edges.

$(S^k)^{(1)}$ is a {\it strictly balanced}
graph, meaning that the ratio of edges to vertices is strictly
smaller for every proper subgraph. A standard result in
random graph theory \cite{Bollo} gives that $n^{-2(k+1)/({2(k+1) \choose
2}-(k+1))}=n^{-1/k}$ is a sharp threshold function for $G(n,p)$
containing a $(S^k)^{(1)}$ subgraph. In particular, if $p^{k}n
\rightarrow \infty$, $G(n,p)$ a.a. contains such a subgraph.

With notation as above, let $S=\{u_1,u_2, \ldots, u_{k+1}\}
\cup \{v_1,v_2, \ldots, v_{k+1} \}$ be the vertices of such a
subgraph. The conditional probability that vertices $\{ u_1,u_2,
\ldots, u_{k+1} \}$ have a common neighbor is no more than
$$ (k+1)p + (n-2k-2)p^{k+1}=  o(1),$$ since $p^{k}n \to 0$
 (so $p \to 0$) and $p^{k+1}n \to 0$.
So a.a. $G(n,p)$ contains a $(S^k)^{(1)}$ subgraph $S$ such that
$\{ u_1,u_2, \ldots, u_{k+1} \}$ has no common neighbor. Note that
in this case $u_i$ is never adjacent to $v_i$ for any choice of
$i$. Then define a retraction of $\XG$ onto $X(S)$ by defining a map $r:G(n,p) \rightarrow S$ on vertices and
extending simplicially. (In particular, the $(S^k)^{1}$ subgraph
is induced.)

For $x \in S$, set $r(x)=x$ and for $x \notin S$, let $i$ be
chosen so that $x$ is not adjacent to $u_i$ and set $r(x)=u_i$.
Such a choice exists for every $x \notin S$ almost always, by
the above. There's no obstruction to extending $r$ simplicially to
a retraction $\tilde{r}:\XG\rightarrow X(S)$, and $X(S)$ is
homeomorphic to $\mathbb{S}^{k}$.

\section{Betti numbers} \label{Betti 1}

First assume that $p^{k+1}n \rightarrow 0$ and $p^k n \rightarrow
\infty$. We wish to prove Theorem \ref{clique_Theorem 5} and in particular to show that a.a. $\beta_k \sim f_k$. For every simplicial complex $\Delta$, we have the Morse inequality \cite{Hatcher}:

$$ -f_{k-1}+f_k-f_{k+1} \le \beta_k \le f_k. $$ The point is that when $p$ is in this interval, $f_k$ is much larger than $f_{k-1}+f_{k+1}$.

By linearity of expectation, we have

\be  \label{emorse} -E[f_{k-1}]+E[f_{k}]-E[f_{k+1}] \le E[\beta_k] \le E[f_k], \ee
and then expanding each term gives 
$$-{n \choose k}p^{{k \choose 2}}+{n \choose k+1}p^{{k+1 \choose 2}}
-{n \choose k+2}p^{{k+2 \choose 2}} \le  E[\beta_k] \le {n \choose
k+1}p^{{k+1 \choose 2}}.$$

%$${n \choose k+1}p^{{k+1 \choose 2}} \left( -\frac{k+1}{n-k}p^{-k}+1 -\frac{n-k-1}{k+2}p^{k+1}\right)  \le  E[\beta_k] \le {n \choose
%k+1}p^{{k+1 \choose 2}}.$$

Since $p^{k+1}n \rightarrow 0$ and $p^k n \rightarrow \infty$, we also have that ${n \choose k}p^{{k \choose 2}} \to 0$ and ${n \choose k+2}p^{{k+2 \choose 2}} \to 0$. Let $Y_k=-f_{k-1}+f_{k}-f_{k+1}$.  Then we have shown so far that \be \label{esim} E[Y_k] \sim
E[\beta_k] \sim E[f_k]. \ee We strengthen this by applying the Second
Moment Method. A standard application of Chebyshev's inequality \cite{Alon}
gives that if $E[X] \rightarrow \infty$ and $\mbox{Var}[X]=o(E[X]^2)$
then a.a. $X \sim E[X]$.  To prove Theorem \ref{clique_Theorem 4} it
suffices to check that $\mbox{Var}[f_k]=o(E[f_k]^2)$ and
$\mbox{Var}[Y_k]=o(E[Y_k]^2)$.

Let $\mu=E[f_k]$ and we have $$\mu^2={n \choose k+1}^2p^{2{k+1
\choose 2}}$$ and $$\mbox{Var}[f_k]=E[f_k^2]-\mu^2.$$ Label the
$(k+1)$-subsets of $[n]$, $1,2,\ldots, {n \choose k+1}$. Let
$A_i$ be the event that subset $i$ spans a $k$-face in
$\XG$, and $A_i \wedge A_j$ the event that both $A_i$ and
$A_j$ occur. Then

\begin{eqnarray*}
E[f_k^2] & = &\sum_{i=1}^{{n \choose k+1}}\sum_{j=1}^{{n \choose k+1}}\mbox{Pr}[A_i \wedge A_j]\\
&=&{n \choose k+1}\sum_{j=1}^{{n \choose k+1}}\mbox{Pr}[A_1 \wedge A_j],\\
\end{eqnarray*} by symmetry. By grouping together $A_j$ by the size
of their intersections with $A_1$ we have

\begin{eqnarray*}
E[f_k^2]&=&{n \choose k+1}\sum_{m=0}^{k+1} {k+1 \choose m}{n-k-1 \choose k+1-m}
p^{2{k+1 \choose 2}-{m \choose 2}}\\
&=&{n \choose k+1}p^{2{k+1 \choose 2}}\sum_{m=0}^{k+1} {k+1 \choose m}{n-k-1 \choose k+1-m}p^{-{m \choose 2}}\\
& \le &\mu^2+{n \choose k+1}p^{2{k+1 \choose 2}}\sum_{m=1}^{k+1} {k+1 \choose m}{n-k-1 \choose k+1-m} p^{ -{m \choose 2} }, \\
\end{eqnarray*}
since ${n-k-1 \choose k+1} \le {n \choose k+1}$ . Then we have that

\begin{eqnarray*}
\frac{E[f_k^2]-\mu^2}{\mu^2}&\le&\frac{\sum_{m=1}^{k+1} {k+1 \choose m}{n-k-1 \choose
k+1-m}p^{-{m \choose 2}}}{{n \choose k+1}}\\
&=&\sum_{m=1}^{k+1} O(n^{-m}p^{-{m \choose 2}})\\
&=&\sum_{m=1}^{k+1} O((n^{-1}p^{-(m-1)/2})^m)\\
&=&o(1),\\
\end{eqnarray*} since $n^{-1}p^{-k-1}=o(1)$ by assumption. We conclude
that a.a. $f_k \sim E[f_k]$. This did not depend on any assumption about $k$, so we also have that a.a. $-f_{k-1} \sim E(-f_{k-1})$ and $-f_{k+1} \sim E(-f_{k+1})$, and adding these three gives that a.a. $Y_k \sim E(Y_k)$.

By equations \ref{emorse} and \ref{esim}, a.a. $\beta_k \sim E(\beta_k)$. The conclusion is that a.a. $\beta_k \sim f_k$, so this completes
the proof of Theorem \ref{clique_Theorem 5}.\\

%\section{Betti numbers, Part 2} \label{Betti 2}

Now we use discrete Morse Theory to prove Theorem \ref{clique_Theorem 6}, that if $p^{k+1}n
\rightarrow \infty$ or $p^{k}n \rightarrow 0$ as $n\rightarrow
\infty$ then $E[\beta_k] /E[f_k] =o(1)$.  For this, a few definitions are
in order. We will write $\sigma < \tau$ if $\sigma$ is a face of $\tau$ of codimension 1.

\begin{definition} A discrete vector field $V$ of a
simplicial complex $\Delta$ is a collection of pairs of faces of $\Delta$
$\{ \alpha <\beta \}$ such that each face is in at most one pair.
\end{definition}

Given a discrete vector field $V$, a {\it closed $V$-path} is a
sequence of faces

$$\alpha_0 < \beta_0 > \alpha_1 < \beta_1 > \ldots < \beta_{n} > \alpha_{n+1}, $$

\noindent such that $\{ \alpha_i < \beta_i \} \in V$ for $i=0,
\ldots, n$ and $\alpha_{n+1}=\alpha_o$. (Note that $\{ \beta_i > \alpha_{i+1} \} \notin V$ since each face is in at most one pair.) We say that $V$ is a {\it discrete gradient vector field} if there are no closed $V$-paths.

Call any simplex not in any pair in $V$ {\it critical}. The main
theorem of discrete Morse Theory is the following \cite{Forman}.

\begin{theorem}[Forman]\label{morse} Suppose $\Delta$ is a simplicial complex with a discrete gradient vector field $V$.
Then $\Delta$ is homotopy equivalent to a CW complex with
one cell of dimension $k$ for each critical $k$-dimensional simplex.
\end{theorem}

First assume that $p^{k+1}n \to \infty$. Since we are assuming the vertex set of $G(n,p)$ is labeled by $[n]$, we can let this induce a total ordering of the vertices. This induces a lexicographic ordering on the faces of $\XG$. For two faces $\sigma$ and $\tau$ of a simplicial complex, we write $\sigma <_{lex} \tau$ if $\sigma$ comes before $\tau$ in the lexicographic ordering. For any set of faces $S$ let $\lexmin(S)$ denote the lexicograhically first element of $S$.

Define a discrete gradient vector field on $\XG$ as follows.

$$V:=\{ \{ \alpha < \beta \} | \mbox{dim}(\alpha)=k \mbox{ and } \beta=\lexmin (\{ b | \alpha < b \mbox{ and } \alpha  <_{lex} b \}) \}$$

 It is clear that no face is in more than one pair, and there are no closed $V$-paths. Let $\sigma:=\{v_1,v_2,\ldots,v_{k+1} \} \subset [n]$, with the vertices listed in increasing order, and set $m:=v_{k+1}$. Then $\sigma$ is a critical $k$-dimensional face of $\XG$ if and only if $ \sigma \in \XG$ and $\sigma \cup \{ x \} \notin \XG$ for every $x >_{lex} m$. These events are independent by independence of edges in $G(n,p)$. So

$$\prob(\sigma \mbox{ is a critical $k$-face})= p^{ {k +1 \choose 2}}(1-p^{k+1})^{n-m}$$
There are ${i-1 \choose k}$ possible choices for $\sigma$ with $v_{k+1}=i$.
Let the number of critical $k$-faces be denoted by $\tilde{f}_k$. We have

\begin{eqnarray*}
 E(\tilde{f}_k) & = & \sum_{i=m+1}^n  {i-1 \choose k} p^{ {k +1 \choose 2}}(1-p^{k+1})^{n-i}\\
 & \le & {n \choose k} p^{ {k +1 \choose 2}} \sum_{i=m+1}^n (1-p^{k+1})^{n-i}\\
 & \le & {n \choose k} p^{ {k +1 \choose 2}} \sum_{i=-\infty}^n (1-p^{k+1})^{n-i}\\
& = & {n \choose k} p^{ {k +1 \choose 2}} \frac{1}{p^{k+1}},\\
\end{eqnarray*}

so

\begin{eqnarray*}
 \frac{E(\tilde{f}_k)}{E(f_k)} & \le & \frac{{n \choose k} p^{ {k +1 \choose 2}} \frac{1}{p^{k+1}}}{{n \choose k+1}p^{{k+1 \choose 2}} }\\
 &=&O\left( \frac{1}{np^{k+1}} \right) \\
 &=&o(1),\\
 \end{eqnarray*}
since $np^{k+1} \to \infty$. By Theorem \ref{morse}, $\XG$ is homotopy equivalent
to a CW complex with at most $\tilde{f}_k$ faces, and by cellular homology, 
$\beta_k \le \tilde{f}_k$ \cite{Hatcher}. So $E(\beta_k) /E(f_k) \to 0$ and this proves the first part of Theorem \ref{clique_Theorem 6}.

Now assume $np^k \to 0$. For each $k$-face $\tau=\{v_1,v_2,\ldots,v_{k+1}\}$ choose $i(\tau) \in \{1, 2, \ldots, k+1 \}$ randomly, uniformly and independently. We'd like to set $$V = \{ \{ \tau-v_{i(\tau)},\tau \} | \tau \in \XG \mbox{ and } \dim(\tau)=k  \},$$ but this might not be a discrete gradient vector field. There are two things that might go wrong. Some $(k-1)$-faces might be in more than one pair, and there might be closed $V$-paths. If we remove one pair from $V$ for each such bad event though, we are left with a proper discrete gradient vector field, with at most one critical cell for each bad event. So we compute the expected number of bad events.

Each bad event contains at least one pair of $k$-faces of $\XG$ meeting in a $(k-1)$-face, either resulting in a $(k-1)$-face being in more than one pair, or a closed $V$-path. Let $d$ denote the number of such pairs, which is also the number of pairs of $K_{k+1}$ subgraphs in $G(n,p)$ which intersect in exactly $k$ vertices. In such a situation there are $k+2$ vertices and at least ${k +2 \choose 2} -1$ edges total, and given a set of $k+2$ vertices there are ${k +2 \choose 2}$ possible choices for a pair of $K_{k+1}$ intersecting in $k$ vertices, so

\begin{eqnarray*}
E(d) & = & {k+2 \choose 2} {n \choose k+2}p^{ {k+2 \choose 2} -1 }\\
& = & {k+2 \choose 2}  {n \choose k+2}p^{\frac{k(k+3)}{2} }.
\end{eqnarray*}

\noindent Then

\begin{eqnarray*}
\frac{E(d)}{E(f_k)} & = & \frac{  {k+2 \choose 2}{n \choose k+2}p^{\frac{k(k+3)}{2} } } { {n \choose k+1} p^{\frac{k(k+1)}{2}} }\\
& = & O(np^k) \\
& = & o(1), \\
\end{eqnarray*}

\noindent since $np^k \to 0$ by assumption. Again, by Theorem \ref{morse} and cellular homology, $\beta_k \le d$, so this completes the proof of Theorem \ref{clique_Theorem 6}.

\section{Random simplicial complexes} \label{rsc}

$\XG$ seems to us a natural probability space of simplicial complexes to study topologically, in part because every simplicial complex is homeomorphic to a clique complex, e.g. by barycentric subdivision \cite{Bjorner}. But of course there are many other possible definitions of random simplicial complexes. 

Linial and Meshulam give a definition for random $2$-complexes $Y(n,p)$ which ``locally'' look like $G(n,p)$, and exhibited a sharp $\Z_2$-homological analogue of Theorem \ref{clique_Theorem ER} \cite{Nati}. This was subsequently generalized to $d$-dimensional complexes and arbitrary fixed finite coefficients $\Z_m$ by Meshulam and Wallach \cite{Wallach}. In \cite{Hoffman}, it is shown that the threshold for vanishing of $\pi_1(Y(n,p))$ is much larger than the Linial-Meshulam-Wallach threshold for $H_1(Y(n,p),\Z_m)$ . The corresponding question for $H_1(Y(n,p),\Z)$ is, as far as we know, still open.

Pippenger and Schleich study a different sort of random $2$-complexes, made by gluing edges of triangles together randomly \cite{Pipp}. Their complexes are pseudomanifolds, and the main motivation is giving quantitative results about fluctuations in the topology of spacetime, in a discrete version of quantum gravity. 

In another article \cite{Kahle}, we study the {\it neighborhood complex} of a random graph $\mathcal{N}[G(n,p)]$. The results are comparable to what we find here: each fixed homology group is roughly unimodal in $p$, and the  nontrivial homology of a random $d$-complex is concentrated in a small number of dimensions. Applications are discussed to topological bounds on chromatic number.

\section{Future directions} \label{open}
\label{open}

Although Theorem \ref{clique_Theorem 2} is technically a generalization of one direction of Theorem \ref{clique_Theorem ER},
it is not clear if it is best possible and we are of the opinion that
it probably is not. We conjecture that Theorem \ref{clique_Theorem 4} is tight instead, and that if $p=n^{\alpha}$ with $\alpha > -1/(k+1)$ then $\XG$ is $k$-connected, almost always.

In a sense, this would be close to determining the homotopy type of $\XG$ when $-1/k < \alpha < -1/(k+1)$. In particular, if one could establish this conjecture, and also show that $\Homology_k(\XG,\Z)$ is torsion free, then standard results in combinatorial homotopy theory \cite{Bjorner} (Theorem 9.18) would imply that if $p=n^{\alpha}$ with $-1/k <\alpha<-1/(k+1)$ then $\XG$ is a.a. homotopy equivalent to a wedge of $k$-dimensional spheres. However,  note that even showing that $\Homology_k(\XG,\Z)$ is free of $m$-torsion for every fixed $m$ would not be good enough, since it is still possible that there is $m$-torsion, with $m$ tending to infinity along with $n$.

Many simplicial complexes arising in combinatorics are homotopy equivalent to a wedge of spheres, and Robin Forman, among others, have wondered if there is any good reason why \cite{Forman}. Combinatorially defined simplicial complexes frequently arise as order complexes of posets, hence automatically are clique complexes, so hope these results are one step toward giving a reasonable answer this question.

\section*{Acknowledgements}

The author thanks his thesis advisors Eric Babson and Chris Hoffman for many illuminating conversations, without which this work would not have been possible. Also, thanks to Tristram Bogart, Anton Dochtermann, Dmitry Kozlov, Nati Linial, and Roy Meshulam, for insight and many helpful comments. The author gratefully acknowledges support by NSA grant \#H98230-05-1-0053, the University of Washington's NSF-VIGRE grant, and MSRI's program in Computational Applications of Algebraic Topology.

\end{document}